\documentclass{llncs}

\usepackage[hyphens]{url}
\usepackage{cite}
\usepackage{graphicx}
\usepackage{framed}
\usepackage{breqn}
\usepackage{amsfonts}
\usepackage{listings}
\usepackage{color}

\begin{document}

\title{Two almost-circles, and two real ones}
\author{Zolt\'an Kov\'acs}
\institute{
The Private University College of Education of the Diocese of Linz\\
Salesianumweg 3, A-4020 Linz, Austria\\
\email{zoltan@geogebra.org}
}

\maketitle              % typeset the title of the contribution

\begin{abstract}
Implicit locus equations in GeoGebra allow the user to do experiments with generalization of the concept
of ellipses, namely with $n$-ellipses. By experimenting we obtain
a geometric object that is very similar to a set of two circles.

\end{abstract}
\section{GeoGebra: a symbolic tool to generalize concepts}
\textit{GeoGebra} \cite{gg} is a well known dynamic geometry software package with millions
of users worldwide. One of its main purposes is to visualize geometric invariants.
Recently GeoGebra has been supporting investigation of geometric constructions
also symbolically by exploiting the strength of the embedded computer algebra system (CAS) \textit{Giac}
\cite{GiacGG-RICAM2013}. A possible use of the embedded CAS is automated reasoning
\cite{gg-art-doc-gh}. In this paper a particular use of the implicit locus
derivation feature \cite{ART-ISSAC2016} is shown, by using the command
\texttt{LocusEquation} with two inputs: a Boolean expression and the sought mover point.
For example, given an arbitrary triangle $ABC$ with sides $a$, $b$ and $c$,
entering \texttt{LocusEquation($a$==$b$,$C$)} results in the perpendicular bisector $d$ of $AB$,
that is, if $C$ is chosen to be an element of $d$, then the condition $a=b$ is satisfied.

Obtaining implicit loci is a new method in GeoGebra to get interesting facts on classic theorems.
These facts have deep connections to algebraic curves which usually describe generalization of
the classic results. Sometimes it is computationally difficult to obtain the curves
quickly enough, but some recent improvements in Giac's elimination algorithm opened
the road to very effectively investigate a large number of geometric constructions \cite{mcs,ACA2017}
including \textit{Holfeld's 35th problem} \cite{mcs,ACA2015}, a generalization of the
\textit{Steiner-Lehmus theorem} \cite{mcs,SteinerLehmus} or \textit{the right triangle altitude theorem}
\cite{ART-ISSAC2016}, among many others.

We need to admit that the possibility to generalize well known theorems is a consequence
of using \textit{unordered geometry} \cite[p.~97]{chou} in the applied tools and theories. In unordered geometry
one cannot designate only the expression of sums of given quantities like the lengths of
a segment,
so the signed quantities will be considered at the same time. (See \cite[p.~59]{chou} for an example
on irreducible problems and indistinguishable cases.) Therefore we obtain a larger
set of points (that is, an \textit{extended locus}) for the resulting algebraic curve as expected. The obtained set may be inconvenient
in some cases (since the output differs from the expected result),
but it can be still fruitful to get some interesting generalizations.

\section{A generalization of the definition of ellipse}

\begin{figure}
\begin{center}
\includegraphics[width=0.7\textwidth]{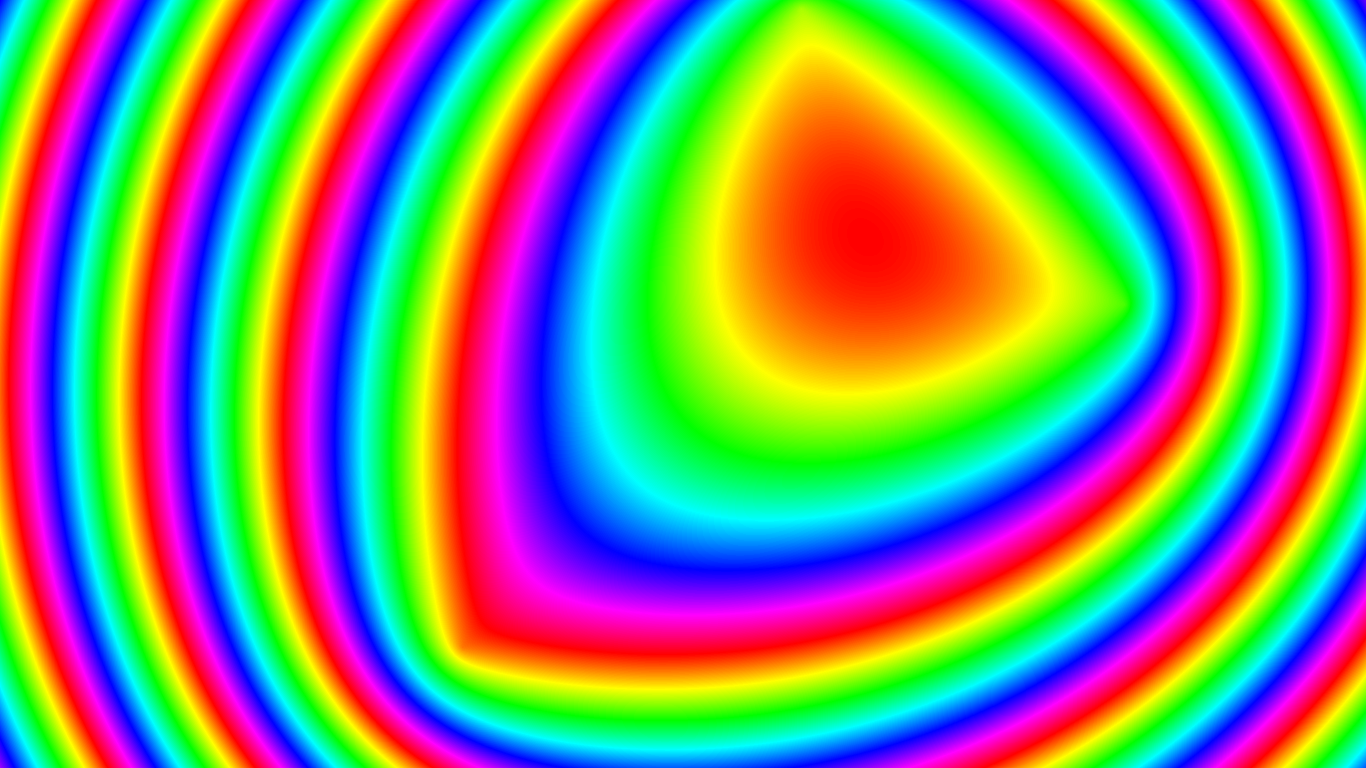}
\caption{A set of curves with fixed points $A_1,A_2,A_3$ and various segments $s$}
\label{Fekete}
\end{center}
\end{figure}

Let us consider fixed points $A_1,A_2,\ldots,A_n$ and a fixed segment with length $s$ in the plane.
We are searching for all points $P$ such that
$$\sum_{i=1}^n |A_iP|=s.$$
In case $n=1$ the obtained set is a circle with center $A_1$ and radius $s$,
in case $n=2$ we obtain an ellipse with foci $A_1$ and $A_2$ and major axis $s$.
In a recent publication \cite{Fekete} by \'Arp\'ad Fekete the geometry of the case $n=3$ is observed by
coloring the various curves for fixed points $A_1,A_2,A_3$ and various segments with length $s$
(see Fig.~\ref{Fekete}\footnote{Recently a CindyJS \cite{cindyjs1,cindyjs2} applet was written by the author that can
produce the same output with just a couple of lines of code, based on a simple statement like
\texttt{colorplot(hue(re(|x-A|+|x-B|+|x-C|)))}. See the examples
\texttt{66_ellipses.html}, \texttt{66_3-ellipses.html} and \texttt{66_4-ellipses.html}
in the folder \texttt{examples/cindygl} at \url{https://github.com/CindyJS/CindyJS}.}).

This idea is, however, already studied by many others. The first appearance is at J.C.~Maxwell's work \cite{maxwell},
but it is called also \textit{$n$-ellipse} by Sekino, and by Nie, Parrilo and Sturmfels \cite{nie-parrilo-sturmfels},
\textit{multifocal ellipse} by Erd\H{o}s and Vincze, \textit{polyellipse}
by Melzak and Forsyth, and \textit{egglipse} by Sahadevan, among others.

\begin{figure}
\begin{center}
\includegraphics[width=0.8\textwidth]{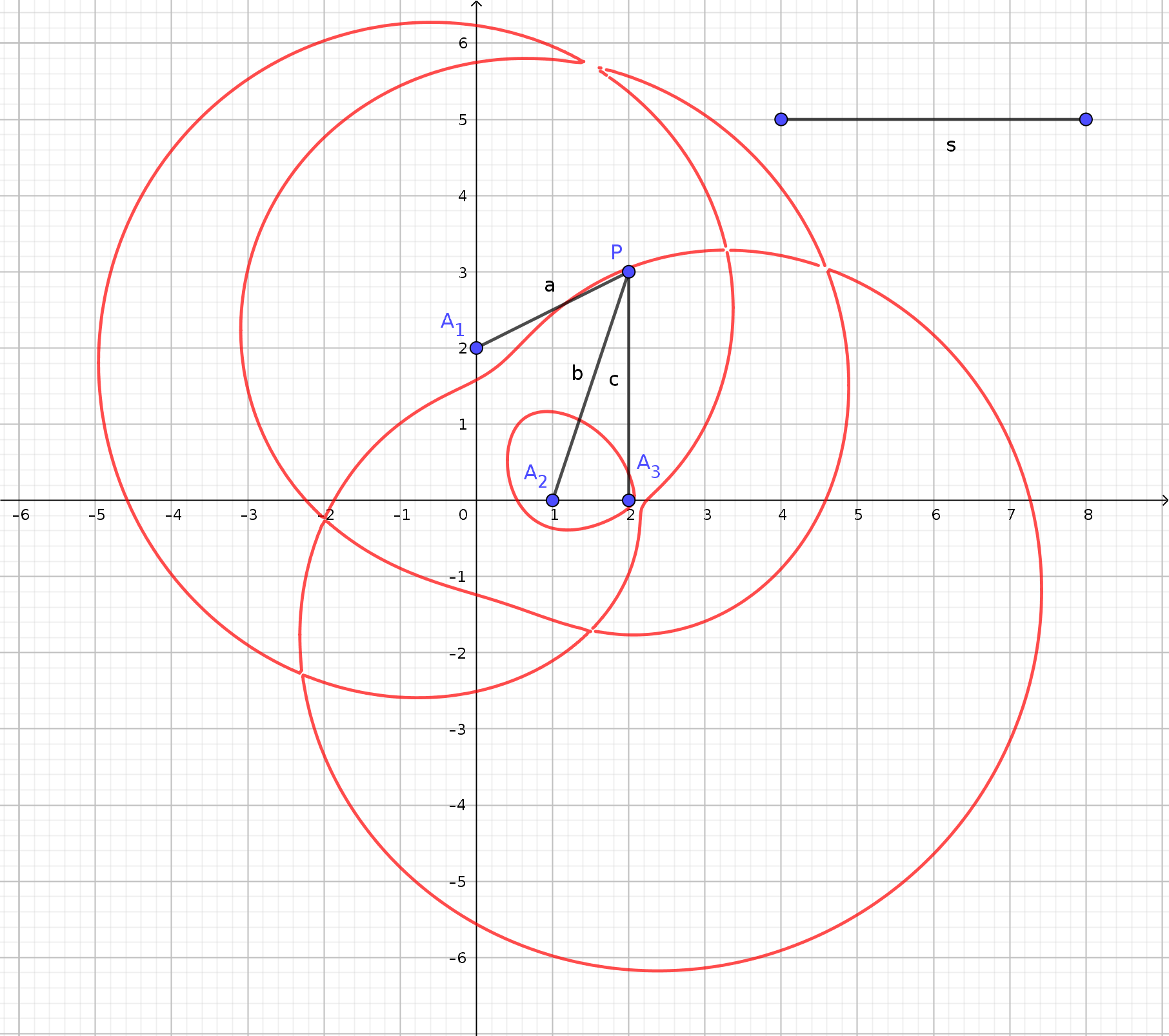}
\caption{GeoGebra's output on \texttt{LocusEquation($a+b+c$==$s$,$P$)}}
\label{ellipse9}
\end{center}
\end{figure}

By using GeoGebra's \texttt{LocusEquation} command one can easily produce a single curve (see Fig.~\ref{ellipse9}),
under the assumption $A_1=(0,2)$, $A_2=(1,0)$, $A_3=(2,0)$ and $s=4$. In this case GeoGebra can also
compute an algebraic equation of the curve, namely
$E(x,y)=9 \; x^{8} + 9 \; y^{8} + 36 \; x^{2} \; y^{6} + 54 \; x^{4} \; y^{4} + 36 \; x^{6} \; y^{2} - 72 \; x^{7} - 48 \; y^{7} - 72 \; x \; y^{6} - 144 \; x^{2} \; y^{5} - 216 \; x^{3} \; y^{4} - 144 \; x^{4} \; y^{3} - 216 \; x^{5} \; y^{2} - 48 \; x^{6} \; y - 220 \; x^{6} - 372 \; y^{6} + 480 \; x \; y^{5} - 964 \; x^{2} \; y^{4} + 960 \; x^{3} \; y^{3} - 812 \; x^{4} \; y^{2} + 480 \; x^{5} \; y + 2136 \; x^{5} + 1712 \; y^{5} + 1656 \; x \; y^{4} + 2080 \; x^{2} \; y^{3} + 3792 \; x^{3} \; y^{2} + 368 \; x^{4} \; y + 446 \; x^{4} + 2846 \; y^{4} - 8256 \; x \; y^{3} + 5452 \; x^{2} \; y^{2} - 7104 \; x^{3} \; y - 14424 \; x^{3} - 6928 \; y^{3} - 22008 \; x \; y^{2} + 688 \; x^{2} \; y + 4980 \; x^{2} + 3132 \; y^{2} + 17376 \; x \; y + 27720 \; x + 3600 \; y - 14175 = 0$.

Visually it is obvious that the output consists of 4 different geometric shapes, so a factorization seems to be useful.
By using \textit{Maple}'s \texttt{evala(AFactor($\ldots$))} command (or 
by using Singular's \cite{Singular} absolute factorization library
\cite{absfact} which is freely available for all readers)
we can learn that the polynomial $E(x,y)$
is irreducible over $\mathbb{C}$. In fact, as \cite[Lemma 2.1]{nie-parrilo-sturmfels} proves,
in all non-degenerate cases the obtained polynomial is of degree $2^3$ and irreducible over $\mathbb{Q}$.

By considering the introductory comments and the paper \cite{mep}, it seems clear that the internal loop belongs to the
expression $a+b+c=s$, while the others to some similar but signed expressions like $a+b-c=s$, $a-b+c=s$ and $-a+b+c=s$.
Theoretically also the expressions $a-b-c=s$, $-a+b-c=s$ and $-a+b-c=s$ could occur, but some geometrical observations
disallow those cases. (In fact, in other setups the latter three cases can also occur, see below in Fig.~\ref{dyncol2}.)
The union of all these curves is our extended locus.

Nie, Parrilo and Sturmfels give  a very similar example as seen in Fig.~\ref{Fekete}
(see \cite[Fig.~3]{nie-parrilo-sturmfels}),
by expressing that all curves are smooth, except those that contain either $A_1$, $A_2$ or $A_3$.
Also, colored set of curves are shown (in 
\cite[Fig.~4]{nie-parrilo-sturmfels}), similarly to Fig.~\ref{ellipse9}.
In addition, they explain why the three extra curves appear in the extended locus. By using the same notions,
we will say that GeoGebra displays the Zariski-closure of a 3-ellipse in Fig.~\ref{ellipse9}.
The length of segment $s$ can also be called \textit{radius} of an $n$-ellipse, and the points $A_1,\ldots,A_n$
will be called its \textit{foci}.

\section{Two almost-circles}

\begin{figure}
\begin{center}
\includegraphics[width=0.5\textwidth]{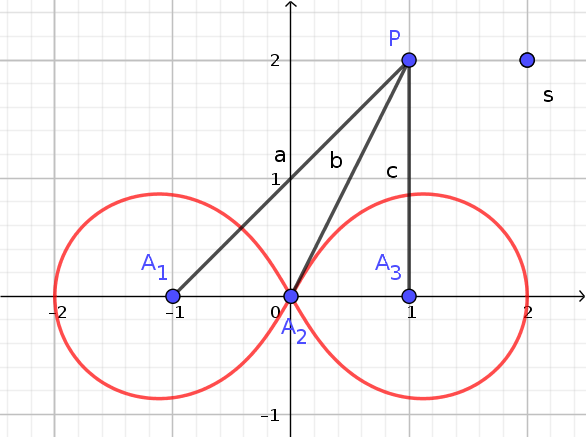}
\caption{An interesting output for $A_1=(-1,0)$, $A_2=(0,0)$, $A_3=(1,0)$, $s=0$, a lemniscate of Booth}
\label{lemniscate}
\end{center}
\end{figure}

By selecting different positions for foci $A_1$, $A_2$, $A_3$ and different radii for $s$ we can obtain a geometrically rich set
of objects. The online GeoGebra applet \url{https://www.geogebra.org/m/tuf3uzf9} can be used by the reader
for own experiments. According to \cite{nie-parrilo-sturmfels}, the  appearing curves are of 8th grade in all cases,
except the degenerate ones, for example, the setup $A_1=(-1,0)$, $A_2=(0,0)$, $A_3=(1,0)$, $s=0$ produces a quartic curve
that looks like a lemniscate (Fig.~\ref{lemniscate}), having the equation $L(x,y)=3x^4+6x^2y^2-12x^2+3y^4+4y^2=0$.
By checking the web page \url{https://en.wikipedia.org/wiki/Hippopede} and learning that the equation $L(x,y)/3=0$
is equivalent to $(x^2+y^2)^2=4x^2-\frac43 y^2$, we identify the obtained curve as a \textit{lemniscate of Booth}.

\begin{figure}
\begin{center}
\includegraphics[width=0.8\textwidth]{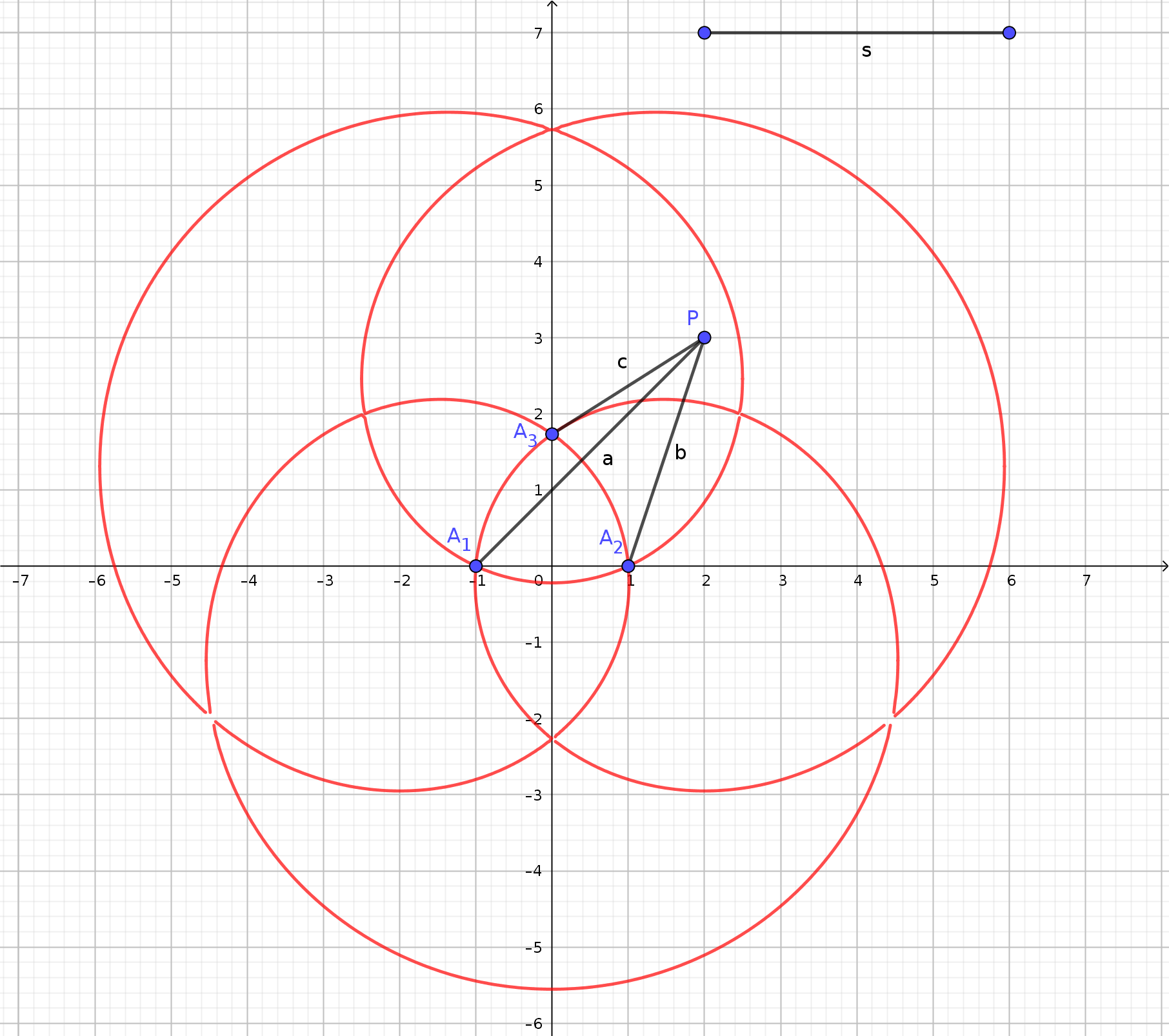}
\caption{An interesting output for $A_1=(-1,0)$, $A_2=(1,0)$, $A_3=(0,\sqrt3)$, $s=4$}
\label{ellipse10}
\end{center}
\end{figure}

\begin{figure}
\begin{center}
\includegraphics[width=0.7\textwidth]{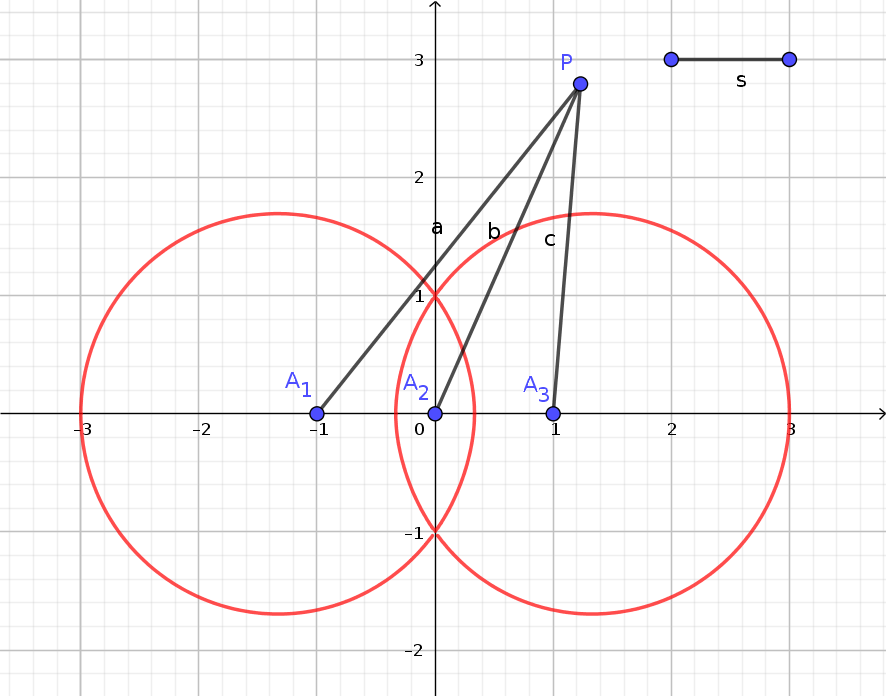}
\caption{Another interesting output for $A_1=(-1,0)$, $A_2=(0,0)$, $A_3=(1,0)$, $s=1$}
\label{ellipse11}
\end{center}
\end{figure}

Some interesting outputs of the 8th grade curve can be found, among others, for
setups $A_1=(-1,0)$, $A_2=(1,0)$, $A_3=(0,\sqrt3)$, $s=4$ (see Fig.~\ref{ellipse10})
and $A_1=(-1,0)$, $A_2=(0,0)$, $A_3=(1,0)$, $s=1$ (see Fig.~\ref{ellipse11}).

\begin{figure}
\begin{center}
\includegraphics[width=0.7\textwidth]{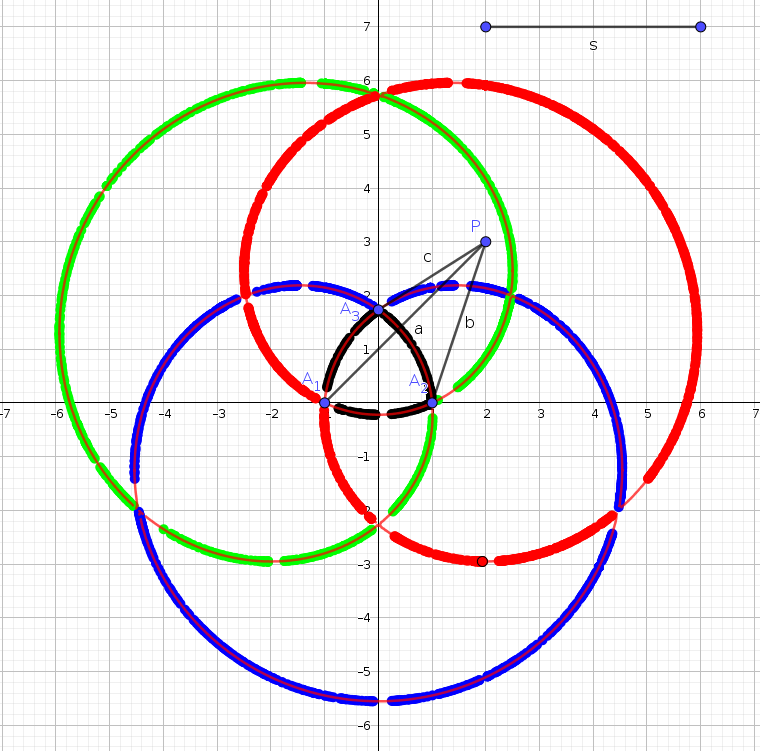}
\caption{Dynamic coloring for $A_1=(-1,0)$, $A_2=(1,0)$, $A_3=(0,\sqrt3)$, $s=4$}
\label{dyncol1}
\end{center}
\end{figure}

\begin{figure}
\begin{center}
\includegraphics[width=0.7\textwidth]{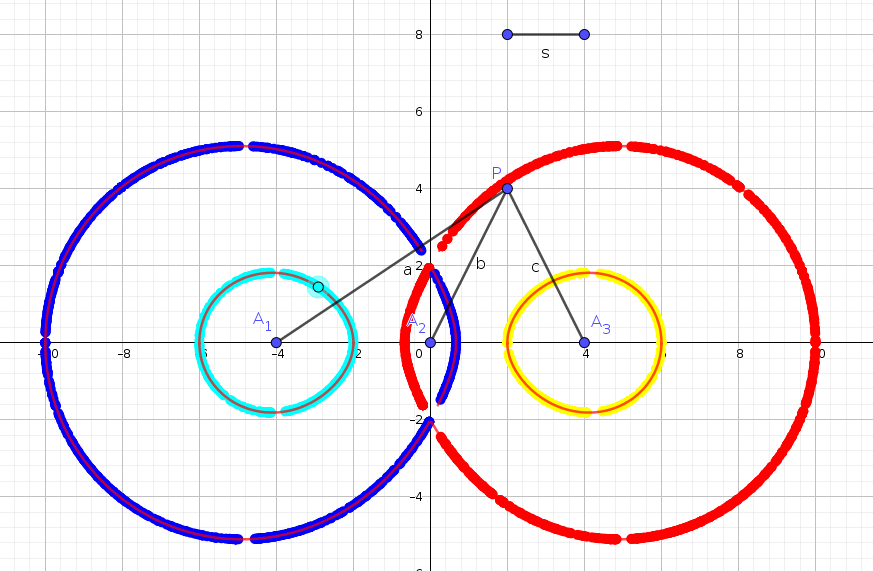}
\caption{Dynamic coloring for $A_1=(-4,0)$, $A_2=(0,0)$, $A_3=(4,0)$, $s=1$}
\label{dyncol2}
\end{center}
\end{figure}

For the former setup (in Fig.~\ref{ellipse10}) we can also learn that, even if the extended curve is just one loop, different parts
of the loop belong to different signed sums. By using a technique ``dynamic coloring'' described in \cite{dyncol}
we can assign RGB components to the signs appearing in the sums, namely, ``red'' for the sign of $|A_1P|$,
``green'' for the sign of $|A_2P|$ and ``blue'' for the sign of $|A_3P|$ in the signed sum $$\sum_{i=1}^3\pm|A_i P|.$$
A $-$ sign adds an RGB component, while a $+$ sign removes it.
For example, the signed sum $+|A_1 P|+|A_2 P|-|A_3P|$ corresponds to $+$red$+$green$-$blue, that is, red and green
are not used, only \textit{blue}. On the other hand, the signed sum $-|A_1 P|-|A_2 P|+|A_3P|$ corresponds to $-$red$-$green$+$blue,
that is, red and green are used (but not blue) which means \textit{yellow} in the RGB system. With this kind of coloring
the black color corresponds to the unsigned 3-ellipse (which is typically the internal loop or the very internal
part of the extended locus). See Fig.~\ref{dyncol1} and \ref{dyncol2}. Also, a GeoGebra applet is available at \url{https://www.geogebra.org/m/bxqa4j9b}
where each setup can be tried out by the user.

Surprisingly enough, the curve in Fig.~\ref{ellipse11} looks like a union of two circles. Factorization, however, does not suggest this conjecture, because
the obtained polynomial, $C(x,y)=
9 \; x^{8} + 9 \; y^{8} + 36 \; x^{2} \; y^{6} + 54 \; x^{4} \; y^{4} + 36 \; x^{6} \; y^{2} - 100 \; x^{6} - 4 \; y^{6} - 108 \; x^{2} \; y^{4} - 204 \; x^{4} \; y^{2} + 182 \; x^{4} - 10 \; y^{4} - 84 \; x^{2} \; y^{2} - 100 \; x^{2} - 4 \; y^{2} +9$,
is irreducible over $\mathbb{C}$.

\begin{figure}
\begin{center}
\includegraphics[width=\textwidth]{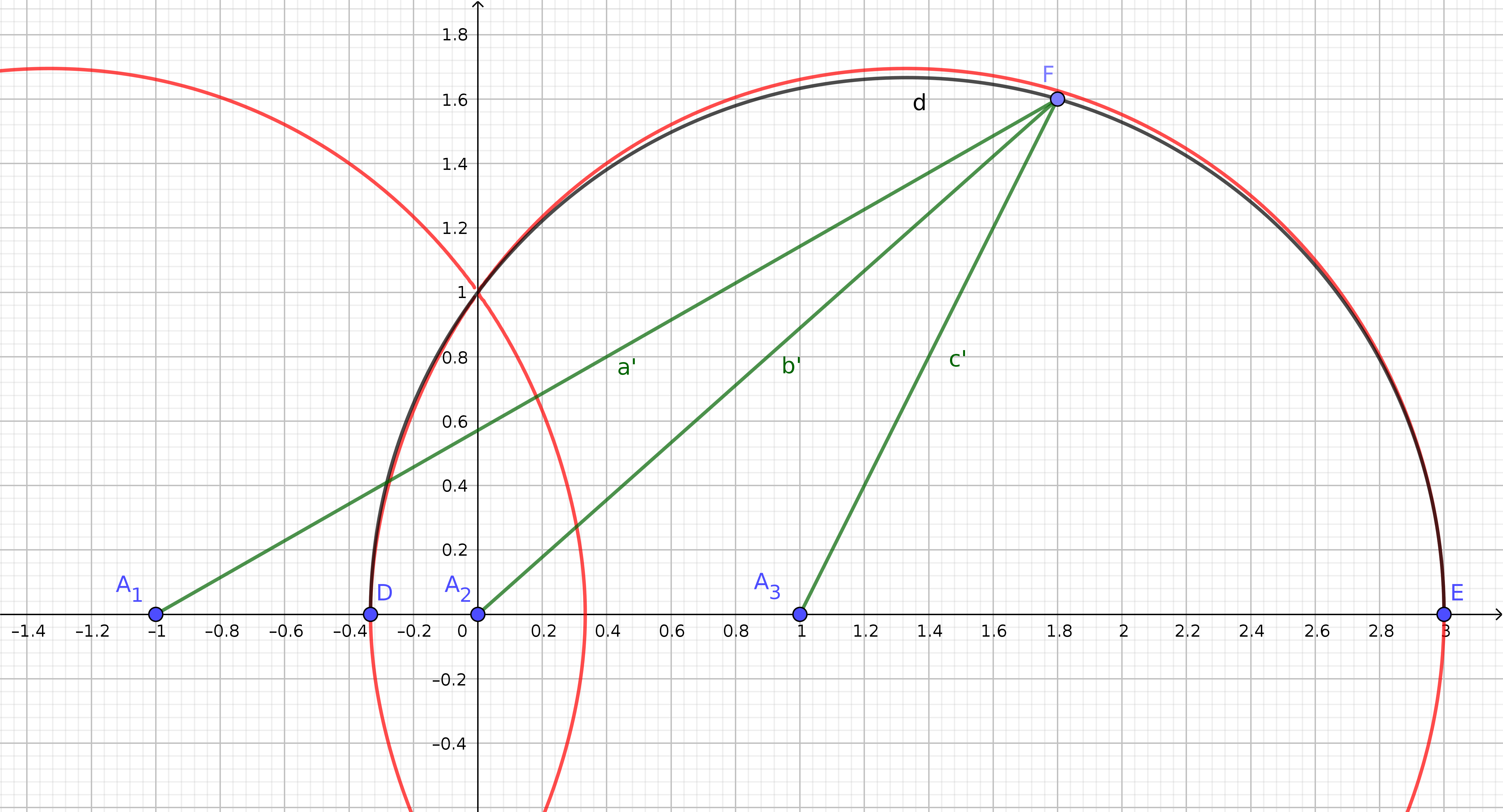}
\caption{A sketch that explains the situation for the two non-circles,
a black semicircle is also drawn to express the difference}
\label{non-circles}
\end{center}
\end{figure}

We can actually prove the fact that Fig.~\ref{ellipse11} does not correspond to two circles. Let us assume the contrary,
that is, the right curve is a circle. Clearly, points $D=(-1/3,0)$ and $E=(1,0)$ are on the extended locus, because
$-|A_1D|+|A_2D|+|A_3D|=-2/3+1/3+4/3=1$ and $-|A_1E|+|A_2E|+|A_3E|=-4+3+2=1$. So we need to assume that the right circle
has center $\left(\frac{3-1/3}2,0\right)=(4/3,0)$ and radius $5/3$. (See Fig.~\ref{non-circles}.)

When considering point $F'=(0,1)$ (this situation is \textit{not} shown in Fig.~\ref{non-circles})
we can still be optimistic. Indeed, $-|A_1F'|+|A_2F'|+|A_3F'|=-\sqrt2+1+\sqrt2=1$, and $F'$ is lying on the circle
since $\sqrt{1^2+(4/3)^2}=5/3$. On the other hand, the point $F=(1.8,1.6)=(9/5,8/5)$ lies on the circle but not on the extended locus:
$\sqrt{\left(\frac{9}{5}-\frac{4}{3}\right)^2+\left(\frac{8}{5}\right)^2}=5/3$, 
$-\sqrt{(1.8-(-1))^2+1.6^2}+\sqrt{1.8^2+1.6^2}+\sqrt{(1.8-1)^2+1.6^2}=-\sqrt{10.4}+\sqrt{5.8}+\sqrt{3.2}\approx0.972270<1$.
Nevertheless, this is still below an error of $3\%$. In general, the error is always below $3.429\%$.

\section{Two real circles}

Viviani's well known theorem for planar triangles states
that \textit{the sum of the distances from any interior point to the sides of an equilateral triangle
equals the length of the triangle's altitude}.

A minor modification of Viviani's theorem can lead to a statement that involves indeed two circles.
We refer here to \cite[Fig.~7]{mcs}
that already reports this result by automated reasoning. Here we repeat a similar kind of proof
by using GeoGebra 5.0.575.0 in three different ways. Namely,
\begin{proposition}\label{prop1}
Let $A_1A_2A_3$ be a regular triangle. The locus of points $P$ such that $|A_1P|+|A_2P|=|A_3P|$ is a
circular arc of the circumcircle of
$A_1A_2A_3$.
\end{proposition}
This proposition is often called \textit{van Schooten's theorem} \cite{van-schooten-wikipedia}.

\begin{figure}
\begin{center}
\fbox{\includegraphics[width=0.8\textwidth]{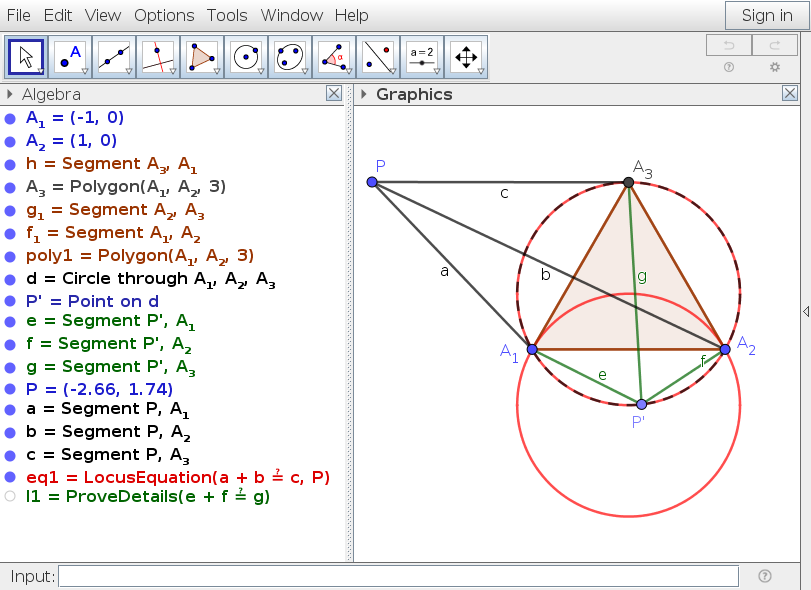}}
\caption{Construction protocol for proving Proposition \ref{prop1} with GeoGebra}
\label{two-circles}
\end{center}
\end{figure}

For a symbolic proof we refer to the GeoGebra construction in Fig.~\ref{two-circles}: a regular triangle is constructed
via the Regular Polygon tool, and the command \texttt{LocusEquation($a+b$==$c$,$P$)} is issued after $P$, $a$, $b$ and $c$
are defined. Notably, the output is a union of \textit{two} circles. The reason behind this is that the regular
triangle is ambiguous: the point $A_3$ can actually be on the other side of segment $A_1A_2$.
The locus equation is $D(x,y)=3 \; x^{4} + 3 \; y^{4} + 6 \; x^{2} \; y^{2} - 6 \; x^{2} - 10 \; y^{2} + 3=0$,
and its factorized form is $$D(x,y)=3\left(x^2+y^2+\frac{2}{3}\sqrt3 y-1\right)\cdot
\left(x^2+y^2-\frac{2}{3}\sqrt3 y-1\right)$$ that clearly corresponds to a union of two circles.
In fact, it is the Zariski-closure of the circular arc over $\mathbb{Q}$.

Another way to prove the same statement in GeoGebra is to draw the circumcircle of a regular triangle $A_1A_2A_3$ and put
a point $P'$ on it, then, by connecting it with $A_1$, $A_2$ and $A_3$, respectively, we get segments with
length $e$, $f$ and $g$. Now the command \texttt{ProveDetails($e+f$==$g$)} gives the output
\texttt{\{true, \{"$A_1=A_2$", "$e=f+g$", "$f=e+g$"\}\}} which can be interpreted in the following complex algebraic geometrical way:

\begin{proposition}\label{prop2}
Let $A_1A_2A_3$ be a regular triangle. The locus of points $P$ such that $|A_1P|+|A_2P|=|A_3P|$ is the circumcircle of
$A_1A_2A_3$, except eventually those points on the circle such that  $A_1=A_2$, $|A_1P|=|A_2P|+|A_3P|$ or $|A_2P|=|A_1P|+|A_3P|$.
\end{proposition}

\begin{figure}
\begin{center}
\fbox{\includegraphics[scale=0.4]{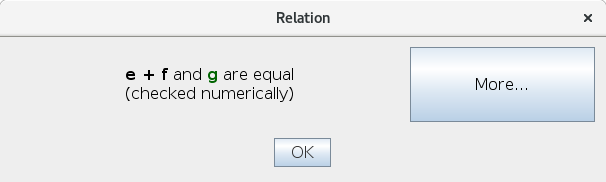}}
\caption{Performing a numerical check in GeoGebra}
\label{rel1}
\end{center}
\end{figure}

\begin{figure}
\begin{center}
\fbox{\includegraphics[scale=0.4]{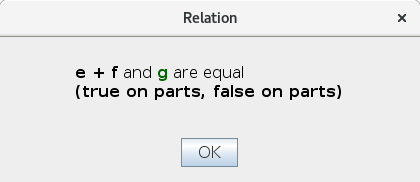}}
\caption{Performing a symbolic check in GeoGebra}
\label{rel2}
\end{center}
\end{figure}

A third method is to enter the command \texttt{Relation($e+f$,$g$)} and obtain the result as shown in Fig.~\ref{rel1},
and then, by pressing the button ``More$\ldots$'', getting the information given in Fig.~\ref{rel2}.

For the meaning of ``true on parts, false on parts'' we refer to \cite{top}, but roughly speaking,
this means that the case $e+f=g$ covers the generally observed equation $e\pm f\pm g=0$ just partly.

At the end of this section we give references to two simple elementary proofs.
One is based on using Ptolemy's theorem and available at Wikipedia \cite{van-schooten-wikipedia}.
The other one is Viglione's idea that was published as a ``proof without words'' in \cite{viglione}.
A GeoGebra applet explaining it in more detailed is available at
\url{https://www.geogebra.org/m/kwgp4abk}.

Finally we note that van Schooten's theorem is about a special case of a $3$-ellipse, namely with
points $A_1=(-1,0)$, $A_2=(1,0)$, $A_3=(0,\sqrt3)$, and radius $s=0$, with the remark that the sum is a \textit{signed} one.

\section{Pedagogical implications}

Today's technology is ready to answer very difficult questions quickly if the problem is entered
in a suitable way. GeoGebra's recent capabilities allow the user---including the student---to find
challenging curves with very little efforts. The output is sometimes surprisingly similar to
well-known geometric objects and the difference cannot be told in a trivial way.

\begin{figure}
\begin{center}
\includegraphics[scale=0.3]{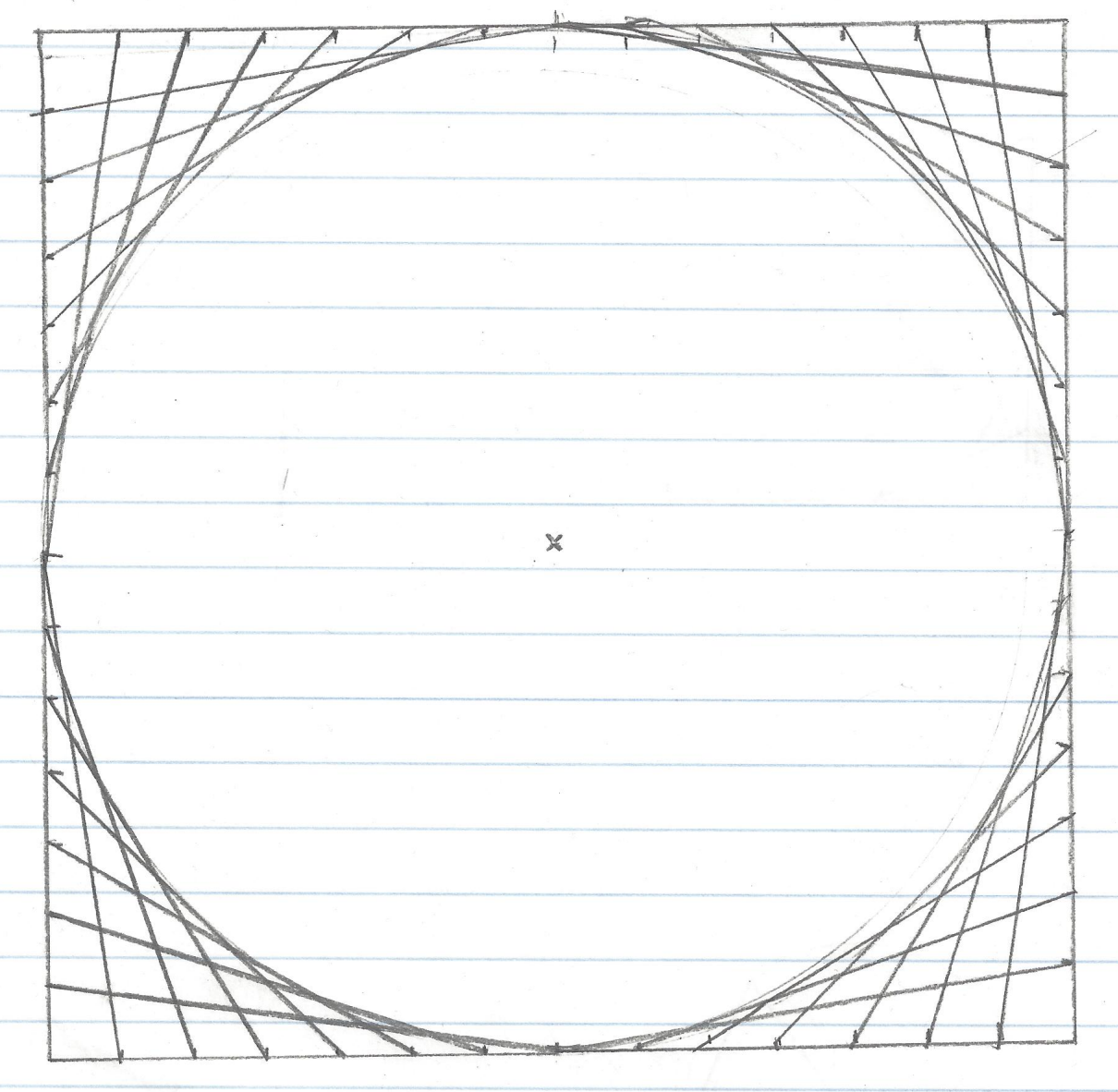}
\caption{Strings that define a union of four parabolas}
\label{4par}
\end{center}
\end{figure}

\begin{figure}
\begin{center}
\includegraphics[scale=0.4]{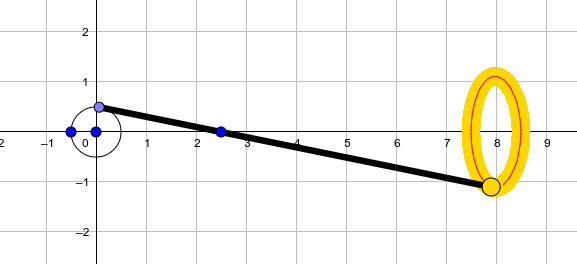}
\caption{A setup of Wittgenstein's rod}
\label{wittg}
\end{center}
\end{figure}

\begin{figure}
\begin{center}
\includegraphics[width=0.6\textwidth]{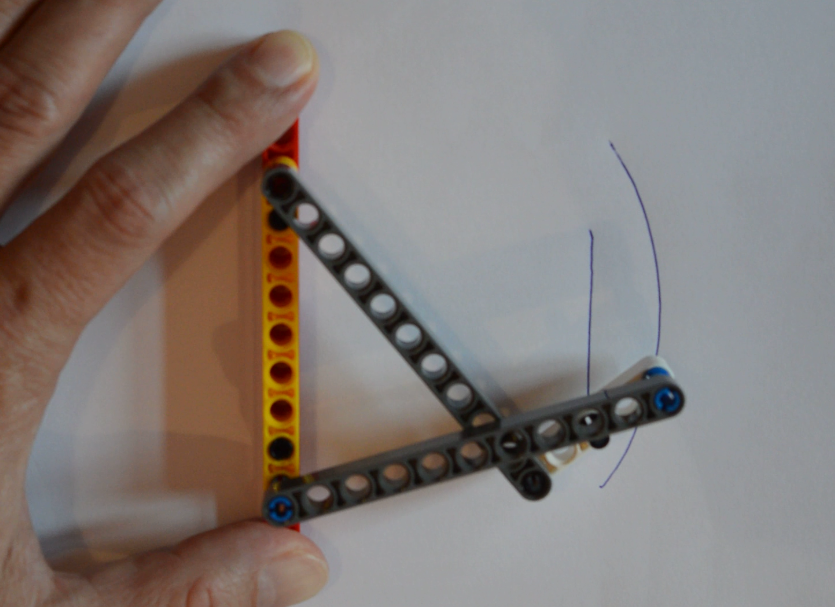}
\caption{Chebyshev's linkage built as a LEGO construction}
\label{chebyshev}
\end{center}
\end{figure}

The same issue may occur on simpler challenges as well. Here we refer to string art parabolas
(that look like a circle, see Fig.~\ref{4par} and \cite{not-a-circle}),
some concrete setups of Wittgenstein's rod (that look like ellipses, but they
are of the 6th degree, see Fig.~\ref{wittg},
available at \url{https://tinyurl.com/wittgensteins-rod}
as a GeoGebra activity) and several other sextic curves that are defined by 4-bar linkages
(they partially look like straight lines, see Fig.~\ref{chebyshev} for an example of both
a seemingly straight line and a seemingly perfect circular arc). For this latter case
we refer to a recent paper \cite{jsc} that describes linkages that can be constructed manually with LEGO parts
and also via a computer program to study their motions.
In fact, all of these almost-curves can be created with just a couple of steps when using a dynamic geometry system
like GeoGebra.

Today's mathematics teachers should be warned about these similarities
and the lack of matches. Luckily, technology is ready enough
to help distinguishing between correct and erroneous conclusions.

\section{Acknowledgments}

The author was partially supported by a grant MTM2017-88796-P from the
Spanish MINECO (Ministerio de Economia y Competitividad) and the ERDF
(European Regional Development Fund).

The author is grateful to \'Arp\'ad Fekete and Noah Dana-Picard for several comments on
a preliminary version of this paper. Special thanks to Carlos Ueno
for his important suggestions that improved the paper significantly.

Finally, many thanks to the CindyJS Team, including Aaron Montag and Michael Strobel, for their support.

\bibliography{kovzol,external}

\end{document}